\documentclass{amsart}
\usepackage{amssymb}
\usepackage{amsthm}
\usepackage{amsmath}
\usepackage{amsfonts}

\newtheorem{theorem}{Theorem}
\newtheorem{corollary}[theorem]{Corollary}

\newtheorem{lemma}[theorem]{Lemma}
\newtheorem{proposition}[theorem]{Proposition}
\newtheorem{definition}[theorem]{Definition}

\def\ppp{{\mathbb{P}}}

\def\qqq{\mathbb{Q}}

\def\ccc{\mathbb{C}}
\def\zzz{\mathbb{Z}}

\def\pf{{\bf Proof}:\ }
\def\qed{$\Box$}
\DeclareMathOperator{\Ind}{Ind}
\DeclareMathOperator{\Res}{Res}

\begin{document}

\title{Representations of finite groups on Riemann-Roch
spaces, II}

\author{David Joyner}
\author{Amy Ksir}
\address{Mathematics Department, United States Naval Academy, Annapolis, MD 21402}
\email{wdj@usna.edu, ksir@usna.edu}
\thanks{The first author was supported in part by an NSA-MSP grant.  The second author was supported in part by a USNA-NARC grant.}

\subjclass[2000]{Primary 14H37}

\date{18 December 2003}

\begin{abstract}
If $G$ is a finite subgroup of the automorphism group of a
projective curve $X$ and $D$ is a divisor on $X$ stabilized by
$G$, then under the assumption that $D$ is nonspecial, we compute
a simplified formula for the trace of the natural representation
of $G$ on Riemann-Roch space $L(D)$.
\end{abstract}

\maketitle

\section{Introduction}

Let $X$ be a smooth projective (irreducible) curve over an
algebraically closed field $k$ and let $G$ be a finite subgroup of
automorphisms of $X$ over $k$. We assume throughout this paper
that either char $k = 0$ or char $k=p$ does not divide the order
of the group $G$.  If $D$ is a divisor of $X$ which $G$ leaves
stable then $G$ acts on the Riemann-Roch space $L(D)$. We are
interested in decomposing this representation into irreducibles.

This question was originally addressed by Hurwitz, in the case
where $D$ was the canonical divisor and $G$ was cyclic, over $k =
\ccc$. Chevalley and Weil expanded this result to any finite $G$
\cite{CW}. Since then further work has been done by Ellingsrud and
L{\o}nsted \cite{EL}, Kani \cite{Ka}, Nakajima \cite{N}, K\"ock
\cite{K1, K2}, and Borne \cite{B}.  In the case where $D$ is a
nonspecial divisor, the character of $L(D)$ has been computed in
the work of Borne \cite{B}. We have computed a simpler formula for
this character, under a rationality criterion.

\begin{theorem}
\label{thm:largedeg} Let $D=\pi^*(D_0)$ be a nonspecial divisor on
$X$ which is a pullback of a divisor $D_0$ on $Y=X/G$ and assume
that the (Brauer) character of $L(D)$ is the character of a
$\qqq[G]$-module. Then for each absolutely irreducible character
of $G$, the multiplicity of the corresponding module $W$ in $L(D)$
is given by

\begin{equation}
\label{eqn:*W} n =
\dim(W)(\deg(D_0)+1-g_Y) - \\
\sum_{\ell =1}^M (\dim(W) - \dim(W^{H_\ell})) \frac{R_\ell }{2}.
\end{equation}
The sum is over all conjugacy classes of cyclic subgroups of $G$,
$H_\ell$ is a representative cyclic subgroup, $W^{H_\ell}$
indicates the dimension of the fixed part of $W$ under the action
of $H_\ell$, and $R_\ell$ denotes the number of branch points in
$Y$ where the inertia group is conjugate to $H_\ell$.
\end{theorem}

One motivation for seeking such a formula comes from coding
theory.  The construction of AG codes uses the Riemann-Roch space
$L(D)$ of a divisor on a curve defined over a finite field.
Automorphisms of $L(D)$ may provide more efficient encoding and
storage of information, for some AG codes.  See \cite{JT} for more
background on AG codes and automorphisms of Riemann-Roch spaces.

In Section 2 we will prove this theorem.  In Section 3 we extend
to the case that $D$ is not necessarily a pullback.  In this case
we use a formula due to Borne \cite{B} which expresses $L(D)$ in
terms of the equivariant degree of $D$ and the ramification module
of the cover, which does not depend on $D$.  Theorem 1 then gives
us a simple formula for the ramification module when it obeys the
rationality condition. This simple formula for the multiplicity of
a $\qqq[G]$-module in the ramification module has also been
obtained by K\"ock \cite{K2} using other methods.  In Section 4,
we give some examples.

The authors would like to thank Bernhard K\"ock for many helpful
discussions while preparing this paper.

\section{Proof of Theorem 1}
\label{sec:2}

We start with some definitions and notation.

Let $X$ be a smooth projective (irreducible) curve over an
algebraically closed field $k$ and let $G$ be a finite subgroup of
automorphisms of $X$ over $k$. We assume that either char $k = 0$
or char $k=p$ does not divide the order of the group $G$. For any
point $P\in X(k)$, let $G_P$ be the inertia group at $P$ (i.e. the
subgroup of $G$ fixing $P$). Our assumptions on char $k$ ensure
that the quotient $\pi :X\rightarrow Y=X/G$ is tamely ramified,
and this group $G_P$ is cyclic.

Let $\langle G\rangle$ denote the set of conjugacy classes of
cyclic subgroups of $G$.  For each class in $\langle G\rangle$
choose a representative cyclic subgroup $H_\ell$, $\ell = 1 \ldots
M$, and partially order them according to the order of the group
so that $H_1$ is the trivial group.   For each branch point of the
cover $\pi: X \to Y$, the inertia groups at the ramification
points $P$ above that branch point will be cyclic and conjugate to
each other. For each $\ell$, let $R_\ell$ denote the number of
branch points in $Y$ where the inertia groups are conjugate to
$H_\ell$.  ($R_1$ may be set to 0; it does not play a role in the
formula).

Let $G_{\qqq}^*$ denote the set of equivalence classes of
irreducible $\qqq[G]$-modules. By results in (\cite{Se}, \S 13.1,
\S 12.4), this set has the same number of elements, $M$, as
$\langle G\rangle$.  For each class in $G_{\qqq}^*$, choose a
representative irreducible $\qqq[G]$-module $V_j$, $j = 1 \ldots
M$, and denote its character by $\chi_j$.  The character table of
$G$ over $\qqq$ is a square matrix with rows labelled by
$G_{\qqq}^*$ and columns labelled by $\langle G\rangle$.  The rows
are linearly independent (as $\qqq$-class functions), so in fact
the character table is an invertible matrix.

Let $F$ be a finite extension of $\qqq$ such that every
irreducible $F[G]$-module is absolutely irreducible (irreducible
over $\ccc$), so that the character table of $G$ over $F$ is the
same as the character table for $G$ over $\ccc$ (\cite{Se}, p.
94).  For each irreducible $\qqq[G]$-module $V_j$, $V_j
\otimes_{\qqq[G]} F[G]$ decomposes into irreducible
$F[G]$-modules. The Galois group of $F$ over $\qqq$ permutes the
components transitively, so each must have the same multiplicity
(the Schur index of the representation $V_j$) and the same
dimension. We write
\begin{equation}
\label{eqn:VWdecomp}
 V_j \otimes_{\qqq[G]} F[G] \simeq m_j
\bigoplus_{r=1}^{d_j} W_{jr},
\end{equation}
where $m_j$ is the Schur index, the $W_{jr}$'s are irreducible
$F[G]$-modules, and $\dim_{\qqq} V_j = m_j d_j \dim_{F} W_{jr}$
for each $r$.  Let $\chi_{jr}$ denote the character of $W_{jr}$.

Theorem \ref{thm:largedeg} is a consequence of the following.

\begin{theorem}
\label{thm:largedegV} Let $D=\pi^*(D_0)$ be a nonspecial divisor
on $X$ and assume that the (Brauer) character of $L(D)$ is the
character of a $\qqq[G]$-module $L(D)_\qqq$. Then for each
irreducible $\qqq[G]$-module $V_j$, its multiplicity in
$L(D)_\qqq$ is given by

\begin{equation}
\label{eqn:*V} n_j = \frac{1}{m_j^2 d_j} \left (
\dim(V_j)(\deg(D_0)+1-g_Y) -  \sum_{\ell =1}^M (\dim(V_j) -
\dim(V_j^{H_\ell})) \frac{R_\ell }{2} \right ).
\end{equation}
\end{theorem}

\pf The proof is similar to the proof of Theorem 2.3 in \cite{Ks}.
 We consider the quotients $X/H_{\ell}$ of $X$ by cyclic subgroups
$H_{\ell}$.  The morphism $\pi:X \to Y$ factors through this
quotient, so on each $X/H_{\ell}$ there is a pullback divisor
$D_{\ell}$ of $D_0$.

First, note that our assumption that $D$ is nonspecial means that
for any quotient $X/H_\ell$, the pullback $D_\ell$ of $D_0$ to
$X/H_\ell$ is also nonspecial. This is because

\[
K_X-D = \pi_\ell ^{*} (K_{X/H_\ell}) + Ram(X/H_\ell) - \pi_\ell
^{*}(D_\ell) = \pi_\ell ^{*}(K_{X/H_\ell}-D_\ell) + Ram(X/H_\ell)
\]
where $Ram(X/H_\ell)$ is the ramification divisor of the covering
$\pi_\ell: X \to X/H_\ell$. Any element of
$L(K_{X/H_\ell}-D_\ell)$ would pull back to $X$ to give an element
of $L(K_X - D - Ram(X/H_\ell))$. Since $Ram(X/H_\ell)$ is
effective, this would also give an element of $L(K_X - D)$,
contradicting our assumption that $D$ is nonspecial.

Now we decompose $L(D)_\qqq$ as

\begin{equation}
\label{eqn:trace} L(D)_\qqq \simeq \bigoplus_{j=1}^M n_j V_j.
\end{equation}
For each $H_\ell$ in $\langle G \rangle$, consider the dimension
of the piece of this module fixed by $H_\ell$. Since the elements
of $L(D)$ fixed by $H_{\ell}$ are exactly the elements of
$L(D_{\ell})$, $\dim_\qqq L(D)_\qqq^{H_\ell} = \dim_k
L(D)^{H_\ell} = \dim_k L(D_\ell)$ and we get an equation for each
$\ell$:

\begin{equation}
\label{eqn:**} \dim_k L(D_\ell) =\sum_{j=1}^M n_j \dim_\qqq (
V_j^{H_\ell }), \ \ \ \ \ 1\leq \ell \leq M.
\end{equation}
This gives us a system of $M$ equations in the $M$ unknowns $n_j$.
We need to show that the matrix $({\rm dim}(V_j^{H_\ell
}))_{j,\ell }$ is invertible, so this system has a unique
solution, and that the above equation is the claimed solution.

First let us consider the matrix $({\rm dim}(V_j^{H_\ell
}))_{j,\ell }$. Each matrix entry is equal to the multiplicity of
the trivial representation of $H_{\ell}$ in the restricted
representation of $H_{\ell}$ on $V_j$. This is the inner product
of characters $\langle \Res_{H_\ell }^G \chi_j, \mathbf{1}\rangle
$, which is defined as

\begin{equation}
\label{eqn:cols} {\rm dim}\, V_j^{H_\ell } = \frac{1}{|H_\ell |}
\sum_{a \in H_\ell } \chi_{j}(a)
\end{equation}

Thus each column of the matrix $({\rm dim}(V_j^{H_\ell }))_{j,\ell
}$ is a sum of columns of the character table of $G$ over $\qqq$.
Each element $a$ in $H_\ell$ generates either all of $H_\ell$ or a
cyclic subgroup of lower order, hence earlier in the list $\langle
G\rangle$. Thus if we write our matrix in terms of the basis of
columns of this character table, we get a lower triangular matrix
with nonzero entries on the diagonal.  This implies that our
matrix is also invertible.

Now it remains to verify that our equation is the correct solution
to (\ref{eqn:**}).

Note that

\begin{equation}
\label{eqn:***} \dim L(D_\ell)= \frac{|G|}{|H|}{\rm
deg}(D_0)+1-g(X/H_{l} ),
\end{equation}
for $1\leq \ell \leq M$, by the Riemann-Roch theorem and the
hypothesis that $D_\ell$ is nonspecial.

We will now substitute (\ref{eqn:*V}) into (\ref{eqn:**}) and
verify that the result agrees with (\ref{eqn:***}), for each
$1\leq \ell \leq M$.  The argument is similar to that in
\cite{Ks}.

Plugging (\ref{eqn:*V}) into (\ref{eqn:**}) gives

{\small{
\[
\begin{split}
\sum_{j=1}^M n_j {\rm dim}(V_j^{H_\ell }) &= ({\rm
deg}(D_0)+1-g_Y)\sum_{j=1}^M \frac{1}{m_j^2 d_j}
{\rm dim}(V_j^{H_\ell }){\rm dim}(V_j) \\
& -\sum_{i=1}^M \left( \sum_{j=1}^M \frac{1}{m_j^2 d_j} [{\rm
dim}(V_j^{H_\ell }){\rm dim}(V_j) - {\rm dim}(V_j^{H_\ell }){\rm
dim}(V_j^{H_i})]
\frac{R_i}{2}\right)\\
\end{split}
\]
}}

Note that
\begin{equation}
\label{eqn:induced}
 {\rm dim}(V_j^{H_\ell }) = \langle
\Res^G_{H_{\ell}} \chi_j, \mathbf{1} \rangle = m_j
\sum_{r=1}^{d_j} \langle \Res^G_{H_\ell} \chi_{jr}, \mathbf{1}
\rangle = m_j \sum_{r=1}^{d_j} \langle \chi_{jr},
\Ind^G_{H_{\ell}} \mathbf{1} \rangle,
\end{equation}
using (\ref{eqn:VWdecomp}) and Frobenius reciprocity. This gives
us

\begin{equation} \begin{split}
\sum_{j=1}^M \frac{1}{m_j^2 d_j} \dim V_j^{H_\ell } \dim V_j &=
\sum_{j=1}^M \frac{\dim V_j}{m_j d_j} \sum_{r=1}^{d_j} \langle
\chi_{jr}, \Ind_{H_\ell}^G \mathbf{1} \rangle \\
&= \sum_{j=1}^M \sum_{r=1}^{d_j} \dim W_{jr} \langle
\Res^G_{H_\ell} \chi_{jr}, \mathbf{1} \rangle \\
&= \frac{1}{|H_\ell|} \sum_{a \in H_\ell} \sum_{j=1}^M
\sum_{r=1}^{d_j} \chi_{jr}(e) \chi_{jr}(a)
\end{split} \end{equation}

The last part of this is summing over all irreducible
$F$-characters of $G$, so the last expression is in fact the inner
product of two columns of the character table for $G$ over $F$.
This inner product will be zero unless $a=e$, so the sum becomes

\begin{equation}
\frac{1}{|H_\ell|} \sum_{j=1}^M \sum_{r=1}^{d_j} \chi_{jr}(e)^2 =
\frac{|G|}{|H_\ell|}.
\end{equation}

We would like to do a similar simplification of

\begin{equation}
\sum_{j=1}^M \frac{1}{m_j^2 d_j} \dim (V_j^{H_\ell }) \dim
(V_j^{H_i})
\end{equation}
using (\ref{eqn:induced}) twice.  The induced representation
$\Ind^G_{H_i} \mathbf{1}$ is the action of $G$ by permutations on
the cosets of $H_{i}$, and thus has a $\qqq[G]$-module structure
as well as an $F[G]$-module structure. It can be decomposed into
irreducible $F[G]$-modules, such that for each $j$ the
multiplicities of the $W_{jr}$'s, $\langle \chi_{jr}, \Ind^G_{H_i}
\mathbf{1} \rangle$, are all equal. Using that fact, Frobenius
reciprocity, and the definition of the Schur inner product, we
have

\begin{equation}
\begin{array}{ll}
\sum_{j=1}^M \frac{1}{m_j^2 d_j} &\dim (V_j^{H_\ell }) \dim
(V_j^{H_i})
\\
&= \sum_{j=1}^M \frac{1}{d_j} \sum_{r=1}^{d_j} \langle
\Res^G_{H_\ell} \chi_{jr}, \mathbf{1} \rangle
\sum_{s=1}^{d_j} \langle \chi_{js}, \Ind_{H_i}^G \mathbf{1} \rangle \\
&= \sum_{j=1}^M \sum_{r=1}^{d_j} \langle \Res^G_{H_\ell}
\chi_{jr}, \mathbf{1} \rangle \langle \chi_{jr}, \Ind_{H_i}^G
\mathbf{1} \rangle \\
&= \sum_{j=1}^M \sum_{r=1}^{d_j} \langle \Res^G_{H_\ell}
\chi_{jr}, \mathbf{1} \rangle
\langle \Res^G_{H_i} \chi_{jr}, \mathbf{1} \rangle\\
&= \frac{1}{|H_\ell|} \frac{1}{|H_i|} \sum_{a \in H_\ell} \sum_{b
\in H_i} \sum_{j=1}^{M} \sum_{r=1}^{d_j} \chi_{jr}(a) \chi_{jr}
(b).
\end{array}
\end{equation}
Again, this last is an inner product of columns of the character
table of $G$ over $k$, so will be zero unless $a$ and $b$ are in
the same conjugacy class. Let $\mathcal{C}_G(a)$ denote the
conjugacy class of $a$ in $G$. We end up with

\begin{equation}
\begin{split}
\sum_{j=1}^M \frac{1}{m_j^2 d_j} \dim (V_j^{H_\ell }) \dim
(V_j^{H_i}) &= \frac{1}{|H_\ell| |H_i|} \sum_{a \in H_\ell}
\#(H_i\cap \mathcal{C}_G(a))
\sum_{j=1}^{M} \sum_{r=1}^{d_j} \chi_{jr}(a)^2 \\
&= |H_\ell\backslash G/H_i|
\end{split} \end{equation}
the number of double cosets.

From this we get

{\small{
\[
\begin{split}
\sum_{j=1}^M n_j \dim V_j^{H_\ell } &=
(\deg(D_0)+1-g_Y)\frac{|G|}{|H_\ell |} -\sum_{i=1}^M
(\frac{|G|}{|H_\ell |}-|H_i\backslash G/H_\ell |)\frac{R_i}{2}\\
&=({\rm deg}(D_0)+1-g_Y)\frac{|G|}{|H_\ell |}+1+
\frac{|G|}{|H_\ell |}(g_Y-1)-g_{X/H_\ell} \\
&={\rm deg}(D_0)\frac{|G|}{|H_\ell |}+1-g_{X/H_\ell}.
\end{split}
\]
}}
where the last equalities come from applying the Hurwitz
formula to the cover $X/H_{\ell} \to Y$ (see \cite{Ks} for
details).   This is (\ref{eqn:***}), as desired. \qed

\textbf{Proof of Theorem \ref{thm:largedeg}.}  We use the
decomposition (\ref{eqn:VWdecomp}) to compute the multiplicity of
each $W_{jr}$ in $L(D)_\qqq \otimes F$.  By our definition of $F$,
each absolutely irreducible character is the character of one of
the $W_{jr}$'s, and the character of $L(D)$ is the same as the
character of $L(D)_\qqq \otimes F$, so this will give us the
correct answer.

The multiplicity of $W_{jr}$ in $V_j$ is $m_j$, and $\dim V_j =
m_j d_j \dim W_{jr}$. Equation \ref{eqn:induced} and the fact that
$\Ind_{H_{\ell}}^{G} \mathbf{1}$ has a $\qqq[G]$-module structure
means that $\dim W_{jr}^{H_\ell}$ is the same for each $r$, so
$\dim V_j^{H_\ell} = m_j d_j \dim W_{jr}^{H_\ell}$. Thus we can
factor $m_j d_j$ out from the inside and multiply the whole thing
by $m_j$ to get formula (\ref{eqn:*W}).  \qed

\textbf{Remark.}  The rationality criterion is necessary for this
formula to be accurate.  If the character of $L(D)$ is not the
character of a $\qqq[G]$-module, it will still be the character of
an $F$-module $L(D)_F$, and $L(D)_F$ will decompose into
irreducibles $W_{jr}$. However in this case for each $j$, the
multiplicities of the $W_{jr}$'s may not be all the same.  The
right hand side of equation (\ref{eqn:*W}) will then compute the
average of these multiplicities:
\begin{small}
\begin{equation}
\label{eqn:LDIrrat}
 \frac{1}{d_j} \sum_{r=1}^{d_j} \langle \chi_{jr}, L(D) \rangle  =
 \dim(W_{jr})(\deg(D_0)+1-g_Y)
  - \sum_{\ell =1}^M
(\dim(W_{jr}) - \dim(W_{jr}^{H_\ell})) \frac{R_\ell }{2}.
\end{equation}
\end{small}

\section{$D$ is not a pullback}

Now we wish to extend our results to the case where $D$ is not
necessarily the pullback of a divisor on $Y=X/G$.  For this we
need to build on work previously done on this problem by Nakajima,
Borne, Ellingsrud and L{\o}nsted, K\"ock, Kani, and others.  We
refer to \cite{B} for references.  We start with two definitions:
the ramification module of the cover $X \to X/G$ and the
equivariant degree of a divisor.

For any point $P\in X(k)$, the inertia group $G_P$ acts on the
cotangent space of $X(k)$ at $P$ by a $k$-character $\psi_P$. This
character is the {\bf ramification character} of $X$ at $P$.  The
ramification module is defined by

\[
\Gamma_G= \sum_{P\in X(k)_{ram}} \Ind_{G_P}^G(\sum_{\ell
=1}^{e_P-1}\ell\psi_P^{\ell}),
\]
where $e_P=|G_P|$. By Theorem 2 in \cite{N}, there is a unique
$G$-module $\tilde{\Gamma}_G$ such that

\[
\Gamma_G= |G| \tilde{\Gamma}_G.
\]
In this paper we are only concerned with $\tilde{\Gamma}_G$, so we
abuse terminology and call $\tilde{\Gamma}_G$ the {\bf
ramification module}.

Now consider a $G$-invariant divisor $D$ on $X(k)$. If
$D=\frac{1}{e_P}\sum_{g\in G}g(P)$ then we call $D$ a {\bf reduced
orbit}. The reduced orbits generate the group of $G$-invariant
divisors $Div(X)^G$.

\begin{definition}
The \textbf{equivariant degree} is a map from $Div(X)^G$ to the
Grothendieck group $R_k(G)=\zzz [G_k^*]$ of virtual $k$-characters
of $G$,

\[
deg_{eq}:Div(X)^G\rightarrow R(G),
\]
defined by the following conditions:

\begin{enumerate}

\item $deg_{eq}$ is additive on $G$-invariant divisors of disjoint
support,

\item If $D=r \frac{1}{e_P} \sum_{g\in G}g(P)$ is an orbit then

\[
deg_{eq}(D)=
\begin{cases}
Ind_{G_P}^G(\sum_{\ell =1}^r\psi_P^{-\ell }), & {\rm if}\ r>0,\\
-Ind_{G_P}^G(\sum_{\ell =0}^{-(r+1)}\psi_P^{\ell }), & {\rm if}\ r<0,\\
0, & r=0,
\end{cases}
\]
where $\psi_P$ is the ramification character of $X$ at $P$.
\end{enumerate}
\end{definition}

\begin{lemma}
\label{lemma:borne} {\bf (Borne's formula)} If $D$ is a
$G$-equivariant nonspecial divisor, then the (virtual) character
of $L(D)$ is given by

\begin{equation}
\label{eqn:borne} \chi(L(D))=(1-g_Y)\chi(k[G])+ deg_{eq}(D) -
\chi(\tilde{\Gamma}_G).
\end{equation}
\end{lemma}

We derive the following from Borne's formula and Theorem
\ref{thm:largedeg}.  The notation is as in Section 1.
\begin{proposition}
\label{prop:rationalrammod} If $\tilde{\Gamma}_G$ has a
$\qqq[G]$-module structure, then it decomposes into irreducible
$\qqq[G]$-modules as

\[
\tilde{\Gamma}_G \simeq \bigoplus_j \frac{1}{m_j^2 d_j} (
\sum_\ell ( \dim(V_j) - \dim(V_j^{H_\ell})) \frac{R_\ell}{2} )
V_j.
\]
\end{proposition}

\pf  The ramification module does not depend on the divisor, so we
compare Theorem \ref{thm:largedeg} with Borne's formula in the
case where $D$ is a pullback.  If $D=\pi^*(D_0)$ is the pull-back
of a divisor $D_0\in Div(Y)$ then the equivariant degree
$deg_{eq}(D)$ has a very simple form. On each orbit, $r$ is a
multiple of $e_P$, so every character of the cyclic group $G_P$
appears.  The equivariant degree on this orbit is induced from a
multiple of the regular representation of $G_P$. Thus we have

\begin{equation}
\label{eqn:eq_deg} deg_{eq}(D)=deg(D_0)\chi(k[G]),
\end{equation}
(This is also a special case of Corollary 3.10 in \cite{B}.)

The first two terms of Borne's formula then become
\begin{equation*}
(\deg D_0 + 1 - g_Y) \chi(k[G]).
\end{equation*}

This is clearly the character of a $\qqq[G]$-module, so $L(D)$
will have a $\qqq[G]$-module structure if and only if
$\tilde{\Gamma}_G$ does.  The rest of the proposition follows from
Theorem \ref{thm:largedeg}. \qed

Proposition \ref{prop:rationalrammod} has also been proven by
K\"ock \cite{K2}, using a different method.

\begin{corollary}
\label{corollary} Suppose that $\tilde{\Gamma}_G$ has a
$\qqq[G]$-module structure.  Let $W$ be an irreducible
$F[G]$-module. Then the multiplicity of the character of $W$ in
$\tilde{\Gamma}_G$ is
\begin{equation}
\sum_\ell ( \dim(W) - \dim(W^{H_\ell})) \frac{R_\ell}{2}.
\end{equation}
\end{corollary}

\pf  The same as the proof of Theorem \ref{thm:largedeg} from
Theorem \ref{thm:largedegV}.  \qed

\textbf{Remark.}  Again, the rationality criterion is necessary.
If $\tilde{\Gamma}_G$ does not have a $\qqq[G]$-module structure,
we get an average of multiplicities, similar to
(\ref{eqn:LDIrrat}):
\begin{equation}
\label{eqn:irratGamma}
\frac{1}{d_j} \sum_{r=1}^{d_j} \langle
\chi_{jr}, \tilde{\Gamma}_G \rangle = \sum_{\ell =1}^M (\dim
(W_{jr}) - \dim (W_{jr}^{H_\ell }))\frac{R_\ell }{2}
\end{equation}
with notation is as in (\ref{eqn:VWdecomp}).

\section{Examples}

\textbf{Example 1.} Consider the nonsingular projective curve X
which is the closure of

\[
 \{ (x,y,t) \in \ccc^3 \ | \  y^2=x(x-2)(x-4),\
t^2=x+4 \}.
\]
This has an action of $G=C_2\times C_2$ given by

\[
\begin{array}{c}
\alpha :(x,y,z)\longmapsto (x,-y,t),\\
\beta :(x,y,z)\longmapsto (x,y,-t),\\
\alpha\beta :(x,y,z)\longmapsto (x,-y,-t).
\end{array}
\]

The quotient by $\beta$ is a degree two cover of an elliptic
curve, ramified at the two points with $x=-4$, so $X$ has genus 2.
The quotient $Y=X/G$ is the projective $x$-line.

The divisor

\[
D=(0,0,2)+(0,0,-2)+(-4,8\sqrt{3},0)+(-4,-8\sqrt{3},0)
\]
is $G$-equivariant, and $2D$ is the pullback of the divisor
$D_0={x=0, x=-4}$ on $Y$.  From the Riemann-Roch theorem we know
that $\dim L(2D) = 7$.

First, let us use Theorem \ref{thm:largedeg} to decompose $L(2D)$
into irreducibles.  The cyclic subgroups of $G$ are the trivial
group, $H_1$ and each of the two-element subgroups generated by
$\alpha$, $\beta$, and $\alpha \beta$.  Let us call the last three
$H_{\alpha}$, $H_{\beta}$, and $H_{\alpha \beta}$.  Each is in its
own conjugacy class.

The cover $X\rightarrow X/G$ has $5$ branch points: three with
inertia group $H_{\alpha}$ (at $x=0,2,4$), one with inertia group
$H_{\beta}$ (at $x=-4$), and one with inertia group
$H_{\alpha\beta}$ (at $x=\infty$). This means

\[
R_{\alpha}=3,\ \ \ R_{\beta}=1,\ \ \ R_{\alpha\beta}=1.
\]

The group $G$ has character table
\[
\begin{array}{c|cccc}
   &   1 &   \alpha &    \beta &    \alpha\beta \\ \hline
\chi_1 &   1 &   1 &    1 &    1 \\
\chi_2 &   1 &   1 &    -1 &    -1 \\
\chi_3 &   1 &   -1 &    1 &    -1 \\
\chi_4 &   1 &   -1 &    -1 &    1
\end{array}
\]
Each irreducible representation is one dimensional, and every
$\ccc[G]$-module is a $\qqq[G]$-module, so $d_j$ and the Schur
index $m_j$ are both 1.  The dimension $\dim (V_j^{H_\ell})$ is 1
if the character of $V_j$ is 1 on the generator and 0 otherwise.
From this we get:

\begin{eqnarray*}
n_1 & = & (2+1-0) - 0 = 3 \\
n_{2} & = & (2+1-0) - \frac{1}{2} (R_{\beta} + R_{\alpha \beta}) =
3 -
1 = 2 \\
n_3 & = & (2+1+0) - \frac{1}{2} (R_{\alpha} + R_{\alpha \beta}) =
3 -
2 = 1 \\
n_4 & = & (2+1+0) - \frac{1}{2} (R_{\alpha} + R_{\beta}) = 3 - 2 =
1.
\end{eqnarray*}

Thus the character of $L(2D)$ is $3 \chi_1 + 2 \chi_2 + \chi_3 +
\chi_4$.

Now let us consider $L(D)$.  The Riemann-Roch theorem tells that
that this will be a three dimensional space.  Since $D$ is not a
pullback from $Y$, we cannot use Theorem \ref{thm:largedeg}.
However, the ramification module does have a $\qqq[G]$-module
structure, so we can use Proposition \ref{prop:rationalrammod}
with Borne's formula.  The calculations above tell us that the
ramification module has character $\chi_2 + 2 \chi_3 + 2 \chi_4.$

Now we need to calculate the equivariant degree of $D$. The
divisor consists of two reduced orbits, the orbit of $(0,0,2)$ and
the orbit of $(-4, 8\sqrt{3}, 0)$.  At the first point the inertia
group is $H_{\alpha}$, and at the second point the inertia group
is $H_{\beta}$.  In both cases the ramification character is the
nontrivial character of $C_2$.  Adding the induced characters of
$G$ gives us $\deg_{eq}(D) = \chi_2 + \chi_3 + 2 \chi_4$.

Adding the pieces of Borne's formula, we get the character of
$L(D)$ to be $\chi_1 + \chi_2 + \chi_4$.  In fact, one can check
that the three functions $\{ 1, \frac{1}{t}, \frac{y}{xt} \}$ form
a basis for $L(D)$, and $G$ acts on the three basis elements by
the three respective characters.  \qed

\textbf{Example 2.}  Let $X = \ppp^1$ and let $G$ be a cyclic
group of prime order $q$.  Let $a$ be a generator of $G$, and let
$a$ act on $X$ by $z \longmapsto \zeta z$, where $\zeta$ is a
primitive $q$th root of unity.  The cyclic subgroups of $G$ are
the trivial group and $G$ itself; the irreducible representations
of $G$ over $\qqq$ are the one-dimensional trivial representation
and a $q-1$ dimensional representation $V$.  Let $\psi$ be the
character of $G$ over $\ccc$ whose value on $a$ is $\zeta$; then
the irreducible characters of $G$ over $\ccc$ are the tensor
powers $\psi, \psi^{\otimes 2}, \ldots, \psi^{\otimes q-1},
\psi^{\otimes q} = \mathbf{1}$. The character of $V$ is $\psi +
\psi^2 + \ldots + \psi^{q-1}$.

The cover $X \to X/G$ is totally ramified at 0 and $\infty$. The
ramification module in this case is a $\qqq[G]$-module, so we can
use either Proposition \ref{prop:rationalrammod} or Corollary
\ref{corollary} to find that
\[
\tilde{\Gamma}_G = \psi + \psi^2 + \ldots \psi^{q-1} = V. \text{
\quad  \qed}
\]

The following example illustrates what can happen when the
rationality condition is not met.

\textbf{Example 3.} Let $X$ be the Klein quartic
\[
\{ \ (x,y,z) \in \ppp^2 \ | \ x^3y+y^3z+z^3x=0 \ \}.
\]
We assume that $k$ contains both cube roots of unity and $7^{th}$
roots of unity; let $\omega$ be a primitive cube root of unity and
$\zeta$ be a primitive seventh root of unity.  Let $G$ be the
group generated by
\begin{eqnarray*}
\sigma:  (x:y:z) & \longmapsto & (y:z:x) \\
\tau: (x:y:z) & \longmapsto & (\zeta x:\zeta^4 y:\zeta^2 z)
\end{eqnarray*}
The group $G$ of automorphisms generated by these two actions is
the semi-direct product $C_3 \rtimes C_7$. (This is not the full
automorphism group of this curve.)  $X$ has genus 2, and the
quotient $Y=X/G$ has genus 0 \cite{E}.

The group $G$ has character table\footnote{This was obtained using
\cite{GAP}. Incidentally, there is only one non-cyclic group of
order $21$, up to isomorphism.}:

\[
\begin{array}{c | ccccc}
  &    e &   \sigma &  \tau & \sigma^{-1} & \tau^{-1} \\ \hline
\chi_1 &   1 &              1 &              1 &              1 &              1 \\
\chi_2 &   1 &         \omega^2 &              1 &           \omega &     1 \\
\chi_3 &   1 &           \omega &     1 &     \omega^2 &   1 \\
\chi_4 &   3 &              0 &  \zeta^3+\zeta^5+\zeta^6 &     0 &  \zeta+\zeta^2+\zeta^4 \\
\chi_5 &   3 &     0 & \zeta+\zeta^2+\zeta^4 & 0&
\zeta^3+\zeta^5+\zeta^6
\end{array}
\]
There are two conjugacy classes of nontrivial cyclic subgroups,
with representatives generated by $\sigma$ and $\tau$. Let $H_3 =
\langle\sigma\rangle$ and $H_7 = \langle \tau \rangle$.   The
irreducible representations over $\qqq$ have characters $\chi_1$,
$\chi_2+\chi_3$, and $\chi_4+\chi_5$. Each has Schur index 1.

The points of $X$ fixed by $H_7$ are $P_1=(1:0:0)$, $P_2=(0:1:0)$,
and $P_3=(0,0,1)$.  These form one orbit under $G$, so $R_7=1$.
 There are seven points in the orbit of $(1:\omega:\omega^2)$ and
seven points in the orbit of $(1:\omega^2:\omega)$, all fixed by
cyclic groups of order 3. Since these form two orbits, we have
$R_3=2$.

We now compute
\begin{equation}
\sum_{\ell =1}^M (\dim (W) - \dim (W^{H_\ell })) \frac{R_\ell}{2},
\end{equation}
as in (\ref{eqn:irratGamma}),  for the irreducible representations
over $\ccc$. We find that
\begin{eqnarray*}
\frac{1}{2} \ \langle \chi_2 + \chi_3, \tilde{\Gamma}_G \rangle &
=
& 1 \\
\frac{1}{2} \ \langle \chi_4 + \chi_5, \tilde{\Gamma}_G \rangle &
= & \frac{7}{2}
\end{eqnarray*}
These give the average multiplicities.  In fact one can compute
directly that $\tilde{\Gamma}_G = \chi_2 + \chi_3 + 3 \chi_4 + 4
\chi_5$.

\end{document}